\documentclass[a4paper, 12pt]{article}

\usepackage{amssymb,amsbsy,amsmath,amsfonts,amssymb,amscd}

\newtheorem{theorem}{Theorem}[section]
\newtheorem{lemma}[theorem]{Lemma}
\newtheorem{proposition}[theorem]{Proposition}
\newtheorem{corollary}[theorem]{Corollary}

\newtheorem{criterion}[theorem]{Criterion}

\newcommand{\ds}{\displaystyle}
\newcommand{\GG}{\mathbb G}
\newcommand{\ofr}{{\mathfrak o}}
\newcommand{\pfr}{{\mathfrak p}}

\newcommand{\lra}{\longrightarrow}
\newcommand{\noi}{\noindent}

\newcommand{\St}{{\mathbf S}{\mathbf t}}

\newcommand{\ZZ}{\mathbb Z}

\newcommand{\Sm}{\mathcal S}
\newcommand{\CC}{\mathbb C}

\newcommand{\HH}{\mathcal H}
\newcommand{\bHH}{\bar \HH}
\newcommand{\VV}{\mathcal V}

\newcommand{\XX}{\mathfrak X}

\newcommand{\Afr}{\mathfrak A}

\newcommand{\Pfr}{\mathfrak P}

\newcommand{\Weyl}{\mathbf W}

\newcommand{\tlambda}{{\tilde \lambda}}

\title{Explicit matrix  coefficients and test vectors  for discrete series representations}
\author{Paul  Broussous}

\date{\today}

\begin{document}
\maketitle

\begin{abstract}
 For the discrete series representations of ${\rm GL}(n)$ over a non-archimedean local field $F$, we define a notion of functions similar to "zonal spherical functions" for the unramified principal series. We prove the existence of such functions in the level $0$ case. As for unramified principal series,  they give rise to explicit matrix matrix coefficients. We deduce a local proof of a known criterion of distinction of discrete series, in the level $0$ case, for the Galois symmetric space ${\rm GL}(n,F)/{\rm GL}(n,F_0 )$, for any unramified quadratic extension $F/F_0$. We also exhibit explicit test vectors when these representations are distinguished. 
 
\end{abstract}

 \tableofcontents

 \bigskip

\section*{Introduction}
In this introduction and throughout the article, $F$ denotes a non-archimedean, non-discrete, locally compact field of any characteristic. 
\smallskip

The aim of this work is  twofold. First we show how to exhibit
explicit matrix coefficients for level $0$ discrete series representations of ${\rm  GL}(n,F)$.  As an application, we obtain  new results for the symmetric space ${\rm GL}(n,F)/{\rm GL}(n,F_0 )$, for an unramified quadratic extension $F/F_0$. We give an entirely local proof of a known criterion of distinction for level $0$ the discrete series representations of ${\rm GL}(n,F)$, as well as explicit test vectors when these representations are distinguished. 
\medskip

 The idea to obtain explicit matrix coefficients is to construct ``zonal
 spherical functions'' similar to classical Satake and Macdonald
 spherical functions attached to unramified principal series
 representations (see [MacD], [Sat], [Cass]). We first recall this
 classical framework.

 Let $\mathbb G$ be a connected reductive
 algebraic group defined over our base field $F$. For simplicity sake,
 let us assume that ${\mathbb G}$ is split over $F$. Let $B$ be a
 Borel subgroup of $G={\mathbb G}(F)$ that we write $B=TU$, for a
 maximal $F$-split torus $T\subset B$ and where $U$ is the unipotent
 radical of $B$. Let $\chi$ be a regular smooth unramified character of
 $T$ and $\pi_\chi$ be the unramified principal series representation
 induced from $\chi$. Let $K$ be a good maximal compact subgroup of
 $G$. It is well known that the Hecke algebra $\HH (G,K)$ is
 commutative and that the fixed vector space $\pi_\chi^K$ is
 $1$-dimensional. In particular $\HH (G,K )$ acts via a character
 $\psi_\chi$ on the line $\pi_\chi^K$.  Let  $c_\chi$ be a
 matrix coefficient of $\pi$ attached to non-zero $K$-invariant vectors in
 $\pi$ and its contragredient. Let ${\bar \HH}(G,K)$ be the
 $\CC$-vector space of bi-$K$-invariant functions on $G$ (of {\it any}
 support). It is not
 an algebra, but nevertheless it is a left $\HH (G,K)$-module for the
 action given by convolution. We let
 $$
 {}_{\psi_\chi} {\bar \HH}(G,K)=\{ f\in {\bar \HH}(G,K)\ ; \varphi
 \star f=\psi_\chi (\varphi ).f, \ \varphi \in \HH (G,K)\} \ .
 $$
 \noi Then the space ${}_{\psi_\chi} {\bar \HH}(G,K)$ is known to be
 $1$-dimensional and the element  $\Phi_\chi$ satisfying
 $\Phi_\chi (1)=1$ is called the {\it zonal spherical function}
 attached to $\chi$. It turns out that $c_\pi$ is
 proportional to $\Phi_\chi$ and that an explicit formula for
 $\Phi_\chi$ is known in terms of the so called {\it Macdonald's
   polynomials}. 
 \medskip

 The same kind of ideas adapt to the case of discrete series representations of $G={\rm GL}(n,F)$.  Let $\Sm (G)$ denote the category of smooth complex representations of $G$. Let $(J,\lambda )$ be a {\it simple 
   type} in $G$ in the sense of Bushnell and Kutzko \cite{BK1}. In
 particular $J$ is an open compact subgroup of $G$ and $\lambda$ is an
 irreducible smooth representation of $J$. Let $\Sm_\lambda (G)$ denote the full subcategory of $G$ whose objects are those representations $(\pi ,\VV )$ that are generated as $G$-modules by their $\lambda$-isotypic component $\VV^\lambda$. We denote by $\HH_\lambda$ (resp. $\bHH_\lambda$) the spherical Hecke algebra (resp. the $\CC$-vector space) of compactly supported functions (resp. functions of any support) $\psi$~: $G\lra {\rm End}_\CC ({\tilde W}_\lambda )$ such that
 $$
 \psi (j_1 gj_2 ) = {\tilde\lambda}(j_1 )\circ \psi (g)\circ {\tilde\lambda} (j_2 ), \ g\in G, j_1 ,j_2\in J\ .
 $$
 \noi Here $W_\lambda$ denotes the space of $\lambda$ and $({\tilde \lambda},{\tilde W}_\lambda )$ denotes the contragredient of $(\lambda ,W_\lambda )$. Note that $\bHH_\lambda$ is naturally a left and right $\HH_\lambda$-module for the action by convolution. A basic result of type theory is that we have an equivalence of categories between $\Sm_\lambda (G)$ and the category $\HH_\lambda$-mod of left $\HH_\lambda$-modules. This equivalence maps a representation $(\pi ,\VV )$ to the module $M_\pi ={\rm Hom}_J (\lambda ,\VV )$. 
 
  Now let $(\pi ,\VV )$ be an irreducible discrete series representation belonging to $\Sm_\lambda (G)$. Then the $\CC$-vector space $M_\pi$ is $1$-dimensional so that $(\pi ,\VV )$ gives rise to a character $\chi_\pi$ of $\HH_\lambda$.  Let us denote by ${}_{\chi_\pi} \bHH_\lambda$ the $\CC$-vector space formed of those functions $\Psi\in\bHH_\lambda$ satisfying
  $$
  \varphi\star \Psi = \chi_\pi (\varphi )\, \Psi , \ \varphi\in \HH_\lambda\ .
  $$
  
  \noi In Proposition \ref{mult1}, we prove that ${\rm dim}_\CC\, {}_{\chi_\pi} \bHH_\lambda\leqslant 1$. Moreover if $\Psi\in {}_{\chi_\pi} \bHH_\lambda$, $w\in W_\lambda$, ${\tilde w}\in {\tilde W}_\lambda$, then the function $c_{w,{\tilde w},\Psi}$~: $G\lra \CC$ defined by 
  $$
  c_{w,{\tilde w},\Psi} (g)=\langle \Psi (g).{\tilde w},w\rangle_\lambda , \ g\in G 
  $$
  is a matrix coefficient of $(\pi ,\VV )$ (Proposition \ref{formula1}).  Now in order to produce an explicit non-zero matrix coefficient of $(\pi ,\VV )$, we have to produce an non-zero element of ${}_{\chi_\pi} \bHH_\lambda$ (and this will prove that ${\rm dim}_\CC \, {}_{\chi_\pi} \bHH_\lambda =1$).  In this goal we use the following recipe. 
  \medskip
  
   Assume first that $J=I$ is an Iwahori subgroup and that $\lambda$ is the trivial character of $I$. In that case the irreducible discrete series representations lying in $\Sm_\lambda$ are twists of the Steinberg representation. As explained in the remark following Lemma \ref{Psi0} a very simple formula is known for the Iwahori spherical matrix coefficient of the Steinberg representation and from this it is easy to produce a non-zero element of ${}_{\chi_\pi} \bHH_\lambda$.
   
    In the general case, it is a central result of \cite{BK1} that for any simple type $(J,\lambda )$ there exists a Hecke algebra isomorphism $\HH_\lambda \sim\HH_{\rm Iw}$, where $\HH_{\rm Iw}$ is the Iwahori-Hecke algebra of ${\rm GL}(m,L)$, for some divisor $m$ of $d$ and some finite field extension $L$ of $F$.  Moreover such an isomorphism is made entirely explicit in \cite{BH} in the case of a level $0$ discrete series representation.  In the level $0$ case, we use this explicit isomorphism to {\it guess} a formula for a non-zero $\Psi_0$ in ${}_{\chi_\pi} \bHH_\lambda$ by copying the corresponding formula for the Steinberg representation of ${\rm GL}(m,L)$ (Theorem \ref{TheFormula}).  
  \medskip
  
    Having constructed our particular non-zero element $\Psi_0\in {}_{\chi_\pi} \bHH_\lambda$ for a level $0$ representation $\pi$, if $w\in W$ and ${\tilde w}\in {\tilde W}$ both are non-zero, the matrix coefficient $c_{w,{\tilde w},\Psi_0}$ is non-zero. More precisely (Proposition \ref{formula1}) there exist $J$-embeddings $W\subset \VV$ and  ${\tilde W}\subset {\tilde \VV}$ such that $c_{w,{\tilde w},\Psi_0}$ is the matrix coefficient of $\pi$ attached to the vector $w\in \VV$ and to the linear form ${\tilde w}\in {\tilde \VV}$. \medskip
    
    In the second part of this work, we apply the previous construction to the theory of distinguished representations. Let us fix a quadratic unramified extension $F/F_0$ and set 
  $G_0 ={\rm GL}(n,F_0 )$. Recall that a representation $(\pi ,\VV )$ of $G$ is said to be $G_0$-distinguished if the intertwining space ${\rm Hom}_{G_0}(\VV ,\CC)$ is non-trivial, where $\CC$ is acted upon by $G_0$ via the trivial character. Let us observe  if $\pi$ is $G_0$-distinguished, its central character $\omega_\pi$ is trivial on the center $F_0^\times$ of $G_0$.  A criterion is known (\cite{An} Theorem 1.3 and \cite{Ma} Corollary 4.2)   to decide whether a discrete series representation of $G$ is $G_0$-distinguished. The proof of this criterion  contains a global ingredient and the aim of this second part is to give a local proof in the case of level $0$ discrete series representations. 
  
   In this aim we use the fact that the symmetric space $G/G_0$ is {\it strongly discrete} in the sense of Sakellaridis and Venkatesh \cite{SV}. This means that if $\pi$  an irreducible discrete series representation of $G$ with central character trivial on $F_{0}^\times$, then any matrix coefficient of $\pi$ is integrable over $G_0 /F_0^\times$. From this it follows that if such a representation $\pi$ has a matrix coefficient $c$ satisfying 
   $$
   \int_{G_0 /F_0^\times} c(g)\, d\mu_{G_0 /F_0^\times} (\dot{g}) \not= 0
   $$
   it has to be $G_0$-distinguished (here $\mu_{G_0 /F_0^\times}$ is a fixed Haar measure on $G_0 /F_0^\times$). 
   \medskip
   
    Let us fix an irreducible level $0$ discrete series representation $(\pi , \VV )$ of $G$ which satisfies the distinction criterion of \cite{An} and \cite{Ma}. Then one proves that there is a natural choice of simple type $(J,\lambda )$ for $\pi$ such that $J$ is ${\rm Gal}(F/F_0 )$-stable and such that $\pi \simeq \pi\circ\sigma$, where $\sigma$ is the generator of ${\rm Gal}(F/F_0 )$. Setting $J_0 =J\cap G_0$, we then prove that $\lambda$ is $J_0$-distinguished with multiplicity one: ${\rm Hom}_{J_0} (\lambda ,\CC )$ is of dimension $1$. Let us abreviate $W_\lambda =W$ and fix a non-zero element $w$ (resp. $\tilde w$) of the line $W^{J_0}$ (resp. of the line ${\tilde W}^{J_0}$). Let $c_0$ be the matrix coefficient $c_{w,{\tilde w}, \Psi_0}$ constructed in the first part of this work. 
    
     We prove the following formula (Proposition \ref{FormulaInt})
     
     $$
     \int_{G_0 /F_0^\times} c_0 (g)\, d\mu_{G_0 /F_0^\times} (\dot{g}) =e\, \langle w,{\tilde w}\rangle_{\lambda} \, P_{\Weyl_0} (\frac{-1}{q_0^f})\ ,
     $$
     \noi where $e$ and $f$ are integers depending on $\pi$ (satisfying $ef=n$), $q_0$ is the size of the residue field of $F_0$,  and $P_{\Weyl_0}$ is the Poincar\'e series of the affine Weyl group $\Weyl_0$ of $G_0$.  Moreover we check that $\langle w,{\tilde w}\rangle\not= 0$ (Lemma \ref{dipendra}) and that $P_{W_0}$ does not vanishes at $-1/q_0^f$ (an explicit formula is known for $P_{W_0}$). As a consequence, this gives a local proof of distinction in our specific case (Theorem \ref{DistGal}). 
     
      Now fix any $J$-embeddings $W\subset \VV$ and ${\tilde W}\subset {\tilde \VV}$. As a byproduct of the previous integral formula we obtain :
     \medskip
     
      $\bullet$ an explicit local period for $\pi$, that is a non-zero element $\Phi$ of ${\rm Hom}_{G_0}(\pi , \CC )$ given by
      $$
      \Phi\ : \ \VV\ni v\mapsto \int_{G_0 /F_0^\times} \langle {\tilde w},v\rangle_\pi\, d\mu_{G_0 /F_0^\times}(\dot{g}) ,  
      $$
      
      $\bullet$ the fact that $w\in \VV$ is a test vector for $\Phi$ : $\Phi (w)\not= 0$. 
  \medskip
  
   A series of observations are needed at this point. 
   \medskip
   
   (1) The "zonal spherical functions" and attached matrix coefficients that we have defined for discrete series representations of ${\rm GL}(n,F)$ may actually be defined for any type $(J,\lambda )$,  in any reductive $p$-adic group $G$,  as soon as the representation $\pi\in \Sm_\lambda (G)$ satisfies ${\rm dim}\, {\rm Hom}_J (\lambda ,\pi )=1$, that is when $\lambda$ occurs in $\pi$ with multiplicity $1$. On the other hand to construct a particular spherical function,  one needs a description of the spherical Hecke algebra $\HH (G,\lambda )$ in terms of an explicit Hecke algebra isomorphism. This description is not always available. For ${\rm GL}(n)$, it is available for simple types of level $0$ thanks to the explicit computations in \cite{BH}. 
   \medskip
   
   (2) Let $\pi$ be an irreducible discrete series representation of ${\rm GL}(n,F)$ and $(J,\lambda )$ be a simple type for $\pi$. Bushnell and Kutzko attach to $\pi$ a certain extended affine Weyl group $\bf W$ (cf. {\S}2). It is a striking fact (Proposition \ref{formula2}) that the explicit matrix coefficient $c_0$ of $\pi$ constructed in this article has its support contained in $J{\mathbf W}J$. 
   \medskip

   (3) The Steinberg representation of ${\rm GL}(n,F)$ is a particular level $0$ discrete series representation. For certain Galois symmetric spaces $G/H$, the $H$-distinction of the Steinberg representation $\St_G$ of $G$ was studied in \cite{BC}, where a conjecture of D. Prasad's is proved. Moreover in \cite{Br} I proved that if $\St_G$ is $H$-distinguished, the value of a well chosen  local period for $\St_G$ at a well chosen Iwahori fixed vector of $\St_G$ is given by a formula similar to that of Proposition \ref{FormulaInt}, that is involving a Poincar\'e series.  
   \medskip

 (4) Let $\pi$ be a generic irreducible smooth representation of ${\rm GL}(n,F)$ which is ${\rm GL}(n,F_0 )$-distinguished. Then Anandavardhanan  and Matringe proved   (\cite{AM} Theorem 1.1) that the {\it essential vector}  of Jacquet,  Piatetski-Shapiro and Shalika is a test vector of $\pi$.  It generally differs from the test vectors constructed in this paper. 
   \bigskip

   The paper is organized as follows. The material needed from "type theory"
is recalled in {\S}1. Level $0$ simple types are given in {\S}2 where one recalls the structure of their Hecke algebras.   The recipe to construct explicit matrix coefficients via "zonal spherical functions" is given in {\S}3. The short section {\S}4 is devoted to the notion of strongly discrete symmetric space and the main result of the article is proved in {\S}5 by integrating the explicit matrix coefficient over $G_0$. 
\bigskip

 The author would like to warmly thank Nadir Matringe and Wiza Silem for usefull discussions, as well as Dipendra Prasad for providing the statment and proof of Lemma \ref{dipendra},  and for helping to improve the containt of this paper.

\section{Review of type theory}

The references for this section are {\S}4 of
\cite{BK1} and {\S}{\S}2-4 of \cite{BK2}.

Let $\GG$ be a connected reductive group
defined over $F$. Fix a Haar measure $\mu_G$ on the locally compact totally
disconnected group $G=\GG (F)$. All representations of $G$ are over complex vector spaces.
We denote by $\Sm (G)$ the abelian  category of smooth complex representations of $G$.

We consider {\it pairs}  $(J,\lambda )$ formed of a
compact open subgroup $J$ of $G$ and an irreducible complex smooth representation
$(\lambda , W)$ of $J$. In particular the representation space $W$ is finite dimensional.
If $(\pi, \VV )$ is an object of
$\Sm (G)$, we denote by $\VV^\lambda$ the $\lambda$-isotypic component of $\VV$. We let
$\Sm_{(J,\lambda )}(G)$ denote the full subcategory of $\Sm (G)$ whose object are those representations
$(\pi ,\VV)$ that are generated by their $\lambda$-isotypic components. Following  Bushnell and Kutzko
\cite{BK2}, one says that $(J,\lambda )$ is a {\it type} if the category $\Sm_{(J,\lambda )}$ is stable
by the operation of taking subquotients.

Recall that if $(K,\rho )$ is a pair in $G$ as above,  one defines two kinds of induction. The smooth
induced representation, denoted ${\rm Ind}_K^G \, \rho$,  acts by right translation on the space
of smooth\footnote{That is smooth under the right action of $G$.}  functions $f$~: $G\lra W_\rho$
satisfying $f(kg)=\rho (k)f(g)$, $g\in G$, $k\in K$, where $W_\rho$ denotes the space of $\rho$.
The compactly induced representation ${\rm ind}_K^G\, \rho$ acts by right translation on
the space of functions $f$~: $G\lra W_\rho$, compactly supported modulo $H$, and satisfying $f(kg)
=\rho (k)f(g)$, $g\in G$, $k\in K$.  Later we shall need the basic fact that the smooth dual
(or contragredient)
of ${\rm ind}_K^G \, \rho$ identifies canonically with ${\rm Ind}_K^G {\tilde \rho}$, where
$({\tilde \rho}, {\tilde W}_\rho )$ denotes the contragredient of $(\rho ,W_\rho )$. Indeed we have a
$G$-invariant non-degenerate pairing
$$
{\rm ind}_K^G \, \rho \times {\rm Ind}_K^G \, {\tilde \rho}, \ \ (f, F)\mapsto \langle f,F\rangle
$$
\noi where
$$
\langle f,F\rangle = \sum_{x\in K\backslash G} \langle f(x) , F(x)\rangle_{\rho}
$$
\noi and $\langle - , - \rangle_\rho$ denotes the natural $K$-invariant pairing on
$W\times {\tilde W}_\rho$. 
\medskip

With the notation as above  fix a type $(J,\lambda )$ in $G$. Recall that the {\it spherical Hecke algebra}
attached to this pair is the space $\HH_\lambda$ of compactly supported functions $f$~:
$G\lra {\rm End}_\CC ({\tilde W})$ satisfying
$$
f(kgl)={\tilde \lambda}(k)\circ f(g)\circ {\tilde \lambda }(l), \ \ k,l\in J, \ \ g\in G,
$$
\noi equipped with the convolution product
$$
( f_1 \star f_2 )(g)=\int_G f_1 (x)f_2 (x^{-1}g)\, d\mu_G (x), \ \ f_1 , f_2 \in \HH_\lambda , \ g\in G .
$$

If $a\in {\rm End}_\CC ({\tilde W})$, we denote by ${}^t a\in {\rm End}_\CC (W)$ its
transpose endomorphism. For $f\in \HH_\lambda$, the fonction ${\tilde f}$~: $g\mapsto {}^t f(g^{-1})$
lies in $\HH_{\tilde \lambda}$. The map $f\mapsto {\tilde f}$, $\HH_\lambda \lra \HH_{\tilde \lambda}$,
is an anti-isomorphism of $\CC$-algebras : $\widetilde{f_1\star f_2}={\tilde f}_1\star {\tilde f}_2$,
$f_1$, $f_2\in \HH_\lambda$.

For $\varphi\in \HH_{\tlambda}$ and $f\in {\rm ind}_J^G \lambda$, the convolution production
$\varphi\star f$, given by
$$
\varphi \star f (g)=\int_{G} \varphi (x) f(x^{-1}g)d\mu_G (x) = \mu_G (J)\,
\sum_{x\in G/K} \varphi (x)f(x^{-1} g), \ \ g\in G \, , 
$$
\noi is well defined and lies in ${\rm ind}_K^G \, \lambda$. It follows that ${\rm ind}_J^G \lambda$
is naturally a $\HH_{\tlambda}$-left module, whence a $\HH_{\lambda}$-right module via
$(f,\varphi )\mapsto f\cdot \varphi = {\tilde \varphi}\star f$,
$f\in {\rm ind}_J^G \lambda$, $\varphi\in \HH_\lambda$. 
\medskip

Let $(\pi ,\VV )$ be an object of $\Sm (G)$. Set $\VV_\lambda := {\rm Hom}_J \, (W, \VV )$. By Frobenius
reciprocity for compact induction from an open compact subgroup, $\VV_\lambda$ canonically identifies
with ${\rm Hom}_G \, ({\rm ind}_J^G\, \lambda , \VV )$. It follows that $\VV_\lambda$
is naturally a left $\HH_\lambda$-module. Hence if $\HH_\lambda {\rm -Mod}$ denotes the category of left
unital $\HH_\lambda$-modules, we have a well defined functor
$$
{\mathbf M}_\lambda\ : \ \Sm (G)\lra \HH_\lambda {\rm -Mod}, \ \ (\pi ,\VV )\mapsto \VV_\lambda
$$

On the other hand we have a functor ${\mathbf V}_\lambda$~:
$\HH_\lambda {\rm -Mod}\lra \Sm (G)$ given on objects by
$$
{\mathbf V}_\lambda (M)={\rm ind}_J^G \, \lambda \otimes_{\HH_\lambda} M
$$
\noi where the action of $G$ on ${\mathbf V}_\lambda (M)$ comes from the (left) action of $G$ on
${\rm ind}_J^G\, \lambda$.

One of the central results of Type Theory is the following result
(cf. \cite{BK2} Theorem (4.3)). 

\begin{theorem} Let $(J,\lambda )$ be a type in $G$. With the notation as above, the functors
  ${\mathbf M}_\lambda$ and ${\mathbf V}_\lambda$ restrict to equivalences of categories
  $\Sm_\lambda (G) \simeq \HH_\lambda {\rm -Mod}$ and are quasi-inverses of each other.
\end{theorem}

\section{Types for the level $0$ discrete series of ${\rm GL}_n$}

In this section $\mathbb G$ is the reductive group ${\rm GL}_n$, for some fixed integer
$n\geqslant 2$,  so that $G={\rm GL}_n (F)$.  Recall that an irreducible smooth representation
$(\pi , \VV )$ of $G$ is a discrete series representation if its central character is unitary and its
matrix coefficients are square integrable modulo the center of $G$. Such a representation is obtained as
follows. Fix a pair of positive integers $(e,f)$ such that $ef=n$ and a unitary irreducible
supercuspidal representation $\tau$ of ${\rm GL}_f (F)$. For $a\in \CC$,  set
$\tau^a := \Vert {\rm det}\Vert^a \tau$, where ${\rm det}$ denotes the determinant map of
${\rm GL}_f (F)$ and $\Vert \Vert$ the normalized  absolute value of $F$. Then
$\tau_M := \tau^{(1-e)/2}\otimes \tau^{(3-e)/2}\otimes \cdots \otimes \tau^{(e-1)/2}$ is a
supercuspidal representation of the standard Levi subgroup ${\rm GL}_f (F)^{e\times}$ of $G$.
Define the {\it generalized Steinberg representation} $\St_e (\tau )$ to be unique irreducible quotient
of the representation of $G$ parabolically induced from $\tau_M$\footnote{Here we consider the
  normalized parabolic induction, taking unitary representations to unitary representations}.
It is a standard fact that
any irreducible discrete series representation is of the form $\St_e (\tau )$ for some pair $(e,f)$ and
some unitary irreducible supercuspidal representation $\tau$ of ${\rm GL}_f (F)$. Moreover $\St_e (\tau )$
is of level $0$ if and only if $\tau$ is. We fix once for all
a pair $(e,f)$ and such a an irreducible level $0$ supercuspidal representation $\tau$.  We denote by
$\omega_\tau$  its central character.

\medskip

The aim of this section is to exhibit a type for $\St_e (\tau )$, that is a pair $(J,\lambda )$ such that
$\St_e (\tau )$ is an object of the subcategory $\Sm_\lambda (G)$. We shall then compute the
$\HH_\lambda$-module ${\mathbf M}_\lambda (\St_e (\tau ))$. We closely follow \cite{BH}. 
\medskip

The irreducible level $0$ supercuspidal representations of ${\rm GL}_f (F)$ are obtained as follows. Let
${\bar \lambda}_0$ be an irreducible cuspidal representation of ${\rm GL}_f (k_F )$. Recall that
we have a canonical isomorphism of groups
$$
{\rm GL}_f (\ofr_F )/ (1+\pfr M_f (\ofr_F ))\simeq {\rm GL}_f (k_f )
$$
\noi so that ${\bar \lambda}_0$ inflates
to an irreducible smooth representation $\lambda_0$ of ${\rm GL}_f (\ofr_F )$.
  The group $U_F \subset {\rm GL}_f (\ofr_F )$ acts as a character $\theta_0$ on the representation space of $\lambda_0$. Fixing an
  extension  $\omega_0$ of $\theta_0$ to the center $F^\times$ of ${\rm GL}_f (F)$ determines an extension
  $\Lambda_0$ of $\lambda_0$ to $F^\times {\rm GL}_f (\ofr_F )$.
  Then ${\rm ind}_{F^\times}^{{\rm GL}_f (F)}\Lambda_0$
   is irreducible and supercuspidal, has level $0$,  and any irreducible level $0$ supercuspidal
    representation is obtained this way. We fix a representation $\lambda_0$ as above so that
    $\tau \simeq {\rm ind}_{F^\times}^{{\rm GL}_f (F)} \Lambda_0$, for some extension $\Lambda_0$ of
    $\lambda_0$. Let us remark that this extension is entirely determined by $\tau$ since one must have
    $\omega_0 =\omega_\tau$.

    Let $\Afr$ be the  hereditary $\ofr_F$-order of ${\rm
      M}_n(F)$ whose elements are the $f\times f$
    block matrices $(a_{ij})_{i,j=1,...,e}$, satisfying
    $a_{ij}\in {\rm M}_f (\ofr_F )$ if $j\geqslant i$, and
    $a_{ij}\in \pfr_F  {\rm M}_f (\ofr_F )$ if $j<i$.
    Its Jacobson radical $\Pfr$ is the set of $f\times f$
    block matrices  $(a_{ij})_{i,j=1,...,e}$, satisfying $a_{ij}\in
    {\rm M}_f (\ofr_F )$ if $j > i$, and
    $a_{ij}\in \pfr_F {\rm M}_f (\ofr_F )$ if $j\leqslant i$.
    It is principal,  generated by the element
    $\Pi$ given in $f\times f$ block form by:
    $$
    \Pi = \left(
    \begin{array}{cccccc}
        0_f          & I_f & 0_f    & 0_f      & \cdots & 0_f      \\
        \vdots       & 0_f & I_f    & 0_f        & \cdots & 0_f       \\        
        \vdots       &     &  \ddots   & \ddots    & \ddots       & \vdots   \\
        \vdots       &     &        &   \ddots  &   \ddots  & 0_f      \\
        0_f          &  \cdots &        &   \cdots  &   0_f  &    I_f  \\ 
        \varpi_F I_f & 0_f & \cdots &      &     \cdots& 0_F
     \end{array}
    \right)
    $$
    \noi where $I_f$ is the $f\times f$ identity matrix. 

    We set $J=\Afr^\times \subset G$ and $J^1 =1+\Pfr\subset J$. Then $J^1$ is a normal subgroup of
    $J$ and $J/J^1 \simeq {\rm GL}_f (k_F )^{e\times}$. We denote by $\lambda$ the inflation of
    ${\bar \lambda}_0^{e\otimes}$ to $J$.

    \begin{proposition} The pair $(J,\lambda )$ is a type in $G$ and $\St_e (\tau )$ is an object of the
      category $\Sm_\lambda (G)$.
    \end{proposition}

    \noi {\it Proof}. The pair $(J,\lambda )$ is a {\it simple type}
    in the sense of $\cite{BK1}$. The fact that simple types are types
    if proved in \cite{BK1} Theorem (8.4.2). The second assertion is
    proved in \cite{BH} {\S}3 Lemma 7.
    \medskip

    We now describe the Hecke algebra $\HH_\lambda$ of $(J,\lambda )$.

    We let  $\Weyl_0 \subset {\rm GL}_e (F)$ be the subgroup of
    monomial matrices with determinant $1$ and entries powers of
    $\varpi_F$. If $w\in \Weyl_0$ we denote by the same symbol $w$ the
    $f\times f$ block monomial matrix of ${\rm GL}_{ef}(F)$ given by
    the Kronecker product $I_f \otimes w$.  In this way we see $\Weyl_0$ as a
    subgroup of $G$. The matrix $\Pi$ normalizes $\Weyl_0$ and we denote by  $\Weyl$ the subgroup
    $\Pi^\ZZ \ltimes \Weyl_0$ of $G$ generated by $\Weyl_0$ and
    $\Pi$. For $i=1,...,e-1$, we denote by $s_i$ the permutation
    matrix in $\Weyl_0$ corresponding to the transposition
    $(i\ i+1)$. Moreover we set $s_0 :=\Pi s_1\Pi^{-1} \in \Weyl$. 
    \medskip

    We denote by $X_0$ the representation space of $\lambda_0$ so that $X:=X_0 \otimes \cdots \otimes X_0$
    is the representation space of $\lambda$. We identify the dual ${\tilde X}$ with
    ${\tilde X}_0^{e\otimes}$. For $i=1,...,e-1$,  we denote by ${\tilde t}_i$ the endomorphism of
    ${\tilde X}$ given by
    $$
    {\tilde t}_i ({\tilde v}_1 \otimes  \cdots \otimes {\tilde v}_i \otimes {\tilde v}_{i+1}\otimes
    \cdots \otimes {\tilde v}_e ) = {\tilde v}_1 \otimes  \cdots \otimes {\tilde v}_{i+1}
    \otimes {\tilde v}_{i}\otimes \cdots \otimes {\tilde v}_e \ , 
    $$
    \noi and we denote by ${\tilde \Gamma}$ the endomorphism of $\tilde X$ defined by
    $$
         {\tilde \Gamma}({\tilde v}_1 \otimes     \cdots \otimes {\tilde v}_e )
         ={\tilde v}_2 \otimes \cdots \otimes {\tilde v}_e \otimes {\tilde v}_1 \ .
    $$

     For $i=1,...,e-1$, there exists a unique element of $\HH_\lambda$ with support $Js_i J$ and
         whose value at $s_i$ is ${\tilde t}_i$.  Similarly there exists
         a unique element $\varphi_\Pi$ of $\HH_\lambda$ with support $J\Pi J = \Pi J=J\Pi$ whose
         value at $\Pi$ is $\Gamma$.
         \medskip

         In order to describe the algebra structure of $\HH_\lambda$, we need to introduce the standard
         affine Hecke algebra $\HH (e,q^f )$ of type ${\tilde A}_e$.
         It has the following presentation. It is
         generated by elements $[s_i ]$, $i=1,...,e-1$, $[\Pi ]$ and $[\Pi^{-1}]$ with the following
         relations:
         \smallskip
         
         (i) $[\Pi ][\Pi^{-1}]= [\Pi^{-1}][\Pi ]=1$

         (ii) $([s_i ]+1)([s_i ]-q^f )=0$

         (iii) $[\Pi ]^2 [s_1 ]=[s_{e-1}] [\Pi ]^2$

         (iv) $[\Pi ] [s_i ]=[s_{i-1}][\Pi ]$, $2\leqslant i \leqslant e-1$

         (v) $[s_i ][s_{i+1}][s_i ] =[s_{i+1}][s_i ][s_{i+1}]$, $1\leqslant i\leqslant e-2$

         (vi) $[s_i ][s_j ]=[s_j ][s_i ]$, $1\leqslant i,j\leqslant e-1$, $\vert i-j\vert \geqslant 2$.
         \smallskip

         We set $[s_0 ]:= [\Pi ][s_1 ][\Pi ]^{-1}$
         \smallskip
         
         Recall that $(\Weyl_0 ,\{ s_0 ,s_2 ,...,s_{e-1}\})$ is a
         Coxeter system of type ${\tilde A}_e$. We denote by
         $l$~: $\Weyl_0 \lra \ZZ_{\geqslant 0}$ its length
         function. We extend it to $\Weyl$ by setting $l(\Pi^a w_0
         )=l(w_0 )$, $a\in {\mathbb Z}$, $w_0\in \Weyl_0$.  If $w=\Pi^k w_1 w_2 \cdots w_l$, $k\in \ZZ$,
         $w_i\in \{ s_1 ,s_2 ,...,s_{e-1}\}$, $i=1,...,l$,  is an element of $\Weyl$, where
         $w_1 w_2 \cdots w_l$ is a reduced expression in $\Weyl_0$, then
         the element $[\Pi ]^k [w_1 ]\cdots [w_l ]$ of $\HH (e,q^f )$
         only depends on $w\in \Weyl$. We denote it $[w]$.
         A basis of $\HH (e,q^f )$ as a $\CC$-vector space is given by $([w])_{w\in \Weyl}$.
         \medskip

         The following classical result allows to perform calculations in $\HH (e, q^f )$.

         \begin{proposition} 
           \label{HeckeCalc} For all $i=0,...,e-1$, for all $w\in \Weyl$, we have
           $$
           [s_i ]\star [w] = \left\{
           \begin{array}{lr}
             [s_i w] & {\rm if}  \ l(s_i w)=l(w)+1\\
             q[s_i w] +(q-1)[w] & {\rm if} \ l(s_i w)=l(w)-1
           \end{array}
           \right.
           $$
           \noi Moreover $l(s_i w)=l(w)+1$ if and only if for all reduced expression $w=w_1 ...w_l$, $l=l(w)$,
           $w_i \in \{s_0 ,...,s_{e-1}\}$, $i=1,...,l$,  we have $w_1 \not= s_i$.

         \end{proposition}

         The following result describes most of the structure of $\HH_\lambda$.

         \begin{theorem} \label{IsoHecke}
           (\cite{BK1} Main Theorem (5.6.6)) a) There exists an isomorphism  of $\CC$-algebras
           $$
           \Upsilon~: \ \ \HH_\lambda \lra \HH (e,q^f )
           $$
           \noi which preserves supports in that
           if $\varphi\in \HH_\lambda$ has support in $JwJ$, for some $w\in \Weyl$, then
           $\Upsilon (\varphi )\in \CC [w]$.

           \noi b) If $\Upsilon$ is such a preserving support isomorphism, then for $i=0,...,e-1$,
           $\Upsilon^{-1}([s_i ])$ does not depends on the choice of $\Upsilon$.

           \noi c) The $\CC$-vector space formed of those $\varphi\in
           \HH_\lambda$ whose support is contained in $\Pi J =J\Pi$ is
           $1$-dimensional. 
         \end{theorem}

         A consequence of the previous result is that, for $w\in \Weyl_0$,  the element
         $\varphi_{[w]}:= \Upsilon^{-1} ([w])$ of $\HH_\lambda$ does not depend on the choice of
         $\Upsilon$. 

         By definition of the $\varphi_{[w]}$, $w\in \Weyl_0$, a straightforward consequence of Proposition
         \ref{HeckeCalc} is the following.
         
         \begin{corollary}  \label{HeckeCalc2}
           For all $i=0,...,e-1$, for all $w\in \Weyl_0$, we have
           $$
           \varphi_{[s_i ]}\star \varphi_{[w]} = \left\{
           \begin{array}{lr}
             \varphi_{[s_i w]} & {\rm if}  \ l(s_i w)=l(w)+1\\
             q\varphi_{[s_i w]} +(q-1)\varphi_{[w]} & {\rm if} \ l(s_i w)=l(w)-1
           \end{array}
           \right.
           $$
         \end{corollary}
     
         The elements $\varphi_{[s_i ]}$, $i=1,...,e-1$, were computed by Bushnell and Henniart.
         
         \begin{proposition} a) (\cite{BH} {\S}2 Proposition 3, relation (2.6.1))
           For $i=1,...,e-1$, the function  $\varphi_{[s_i ]}$ has support
           $Js_i J$ and its value at $s_i$ is $\ds \omega_\tau
           (-1)q^{-\frac{f(f-1)}{2}}\, {\tilde t}_i$.

           \noi b) Define an endomorphism ${\tilde t}_0$ of ${\tilde
             X}$ by
           $$
                {\tilde t}_0 ({\tilde v}_1 \otimes \cdots \otimes
                {\tilde v}_e ) = {\tilde v}_e \otimes {\tilde v}_2
                \otimes \cdots {\tilde v}_{e-1}\otimes {\tilde v}_1
                $$
         \noi Then $\varphi_{[s_0 ]}$ has support $Js_0 J$ ans its value
         at $s_0$ is $\ds \omega_\tau (-1)q^{-\frac{f(f-1)}{2}}\, {\tilde t}_0$
         \end{proposition}
 
 \noi {\it Proof}. Let $\Upsilon$~: $\HH_\lambda \lra \HH (e, q^f )$
 be any support preserving isomorphism. By definition we have
 $\varphi_{[s_0 ]} = \Upsilon^{-1} ([\Pi][s_1 ][\Pi ]^{-1}) =
 \Upsilon^{-1}([\Pi ]) \star \varphi_{[s_1 ]}  \star  \Upsilon^{-1}([\Pi
 ])^{-1}$. By Theorem \ref{IsoHecke}(c), $\Upsilon^{-1}([\Pi ])$ is
 proportional to $\varphi_\Pi$, so that $\varphi_{[s_0 ]} =
 \varphi_{\Pi}\star \varphi_{[s_1 ]}\star \varphi_{[\pi ]}^{-1}$. This latter
   function has support $J \Pi Js_1 J\Pi^{-1}J =Js_{0}J$ and its value
   at $s_0$ is indeed  $\ds \omega_{\tau}(-1) q^{\frac{f(f-1)}{2}}
   {\tilde \Gamma}{\tilde t}_1 {\tilde \Gamma}^{-1}$ $=$ $\ds
   \omega_{\tau}(-1) q^{\frac{f(f-1)}{2}} {\tilde t}_0$.
   \medskip

         We shall need later the expression of the images ${\tilde \varphi}_{[s_i ]}$, $i=0,...,e-1$,
         ${\tilde \varphi}_\Pi$, of $\varphi_{[s_i ]}$, $i=0,...,e-1$, and $\varphi_\Pi$ under the
         anti-isomorphism of algebras $\HH_\lambda \lra \HH_{\tilde \lambda}$, $\varphi\mapsto {\tilde \varphi}$.

         \begin{lemma} \label{valeurshecke} For $i=1,...,e-1$,
           define an element $t_i \in {\rm End}_\CC (X)$ by $t_i (v_1 \otimes
           \cdots \otimes v_e ) = v_1 \otimes \cdots \otimes
           v_{i+1}\otimes v_i \otimes \cdots \otimes v_e$. Moreover
           define $t_0$ and $\Gamma\in{\rm End}_\CC (X)$ by 
            $t_0 (v_1 \otimes \cdots v_e )$ $=$ $v_e \otimes v_2
           \otimes \cdots \otimes v_{e-1}\otimes v_1$ and
           $\Gamma (v_1 \otimes \cdots \otimes v_e )$ $=$ $v_e \otimes v_1 \otimes \cdots v_{e-1}$.
           Then for $i=0,...,e-1$, ${\tilde \varphi}_{[s_i ]}$ is the unique element of $\HH_{\tilde \lambda}$ with
           support $Js_i J$ and taking value $\ds \omega_\tau (-1 ) q^{-\frac{f(f-1)}{2}} \, t_i$ at $s_i$. Similarly
           ${\tilde \varphi}_\Pi$ is the unique element of $\HH_{\tilde \lambda}$ with support $\Pi^{-1}J =J\Pi^{-1}$
           and taking value $\Gamma$ at $\Pi^{-1}$.
         \end{lemma}
           
          \noi {\it Proof}. Straightforward calculations. 
           \medskip

         We now describe the $\HH_\lambda$-module structure of $M:= {\mathbf M}_\lambda (\St_e (\tau ))$.

         \begin{theorem} \label{caractere} (\cite{BH} {\S}3
           Proposition 4 and Proposition 6.(2).)  a) The module $M$ is
           $1$-dimensional so
           that $\HH_\lambda$ acts on $M$ via
           a character $\chi$.

           \noi b) The character $\chi$ is given on a set of generators by the relations:
           \smallskip

           (i) $\chi (\varphi_{[s_i ]})=-1$, $i=1,...,e-1$,

           (ii) $\chi (\varphi_{\Pi})=(-1)^{e-1}\, \omega_\tau (\, (-1)^{e-1} \varpi_F\,  )$
         \end{theorem}

         \noi {\it Remarks}. a) The expression giving  $\chi (\varphi_\Pi )$ depends on the choice of
         the uniformizer. This is not a contradiction since $\Pi$ itself depends on the choice of
         $\varpi_F$.

         \noi b)  If $w\in \Weyl_0$ has reduced expression $w_1 w_2 \cdots w_l$, $w_i \in \{ s_0 ,...,s_{e-1}\}$.
         Then $\chi (\varphi_{[w]})=(-1)^l$. In other words, we have the formula
         $$
         \chi (\varphi_{[w]}) =(-1)^{l(w)}, \ w\in \Weyl_0
         $$
         \noi where we recall that $l$ is the length function of the Coxeter system $(\Weyl_0 , \{ s_0 ,...,s_{e-1}\})$. 
         \medskip

         \section{Explicit matrix coefficients}

         In the previous section we exhibited a type $(J,\lambda )$ for our fixed discrete series
         representation $(\pi ,\VV )=\St_e (\tau )$, of $G={\rm GL}_n (F)$. The left $\HH_\lambda$-module
          is $1$-dimensional: we may identify it with the line
         $\CC =\CC_\chi$, $\HH_\lambda$ acting via the character $\chi$ that we described in Theorem \ref{caractere}.

         Recall that we have an isomorphism of $G$-modules
         $$
         \VV \simeq {\rm ind}_J^G\, \lambda \otimes_{\HH_\lambda} \CC_\chi
         $$

         \noi In other words $\VV$ identifies with the quotient of ${\mathfrak X}:= {\rm ind}_J^G\, \lambda$
         by the sub-$G$-module
         $$
         {\mathfrak X}_\chi := {\rm Span}_\CC\, \{\,  {\tilde \varphi}\star f -\chi (\varphi )f\ ;
         \ f\in \XX , \ \varphi\in \HH_\lambda\, \}
         $$

         \noi Hence the contragredient ${\tilde \VV}$ of $\VV$ is given by
         $$
              {\tilde \VV}\simeq \{ \, \Lambda \in {\rm Ind}_J^G \, {\tilde \lambda}\ ; \
              \Lambda ({\tilde \varphi}\star f)=\chi(\varphi )
              \Lambda (f), \ f\in \XX ,\ \varphi\in
              \HH_\lambda\, \} .
              $$

              Let $\varphi\in \HH_\lambda$, $\Lambda\in {\rm Ind}_J^G \, {\tilde \lambda}$
                and $f\in {\rm ind}_J^G\, \lambda$. We have

        \begin{align*}
        \Lambda ({\tilde \varphi}\star f)  & = & \sum_{x\in J\backslash G}
         \sum_{y\in J\backslash G}\langle {\tilde \varphi} (xy^{-1}f(y) , \Lambda (x)\rangle_\lambda \\
         & =  &   \sum_{x\in J\backslash G}\sum_{y\in J\backslash G}\langle f(y) , \varphi (yx^{-1})\Lambda (x)
         \rangle_{\tilde \lambda}\\
         & = & \sum_{y\in J\backslash G} \langle f(y),  \sum_{x\in J\backslash G}\varphi (yx^{-1})\Lambda (x)
         \rangle_{\tilde \lambda}\\
         & = & (\varphi\star \Lambda )(f)
      \end{align*}

        \noi where the convolution $\varphi \star \Lambda$ is defined by the {\it finite} sum:
        $$
        (\varphi \star\Lambda )(g)=\sum_{x\in J\backslash G} \varphi (gx^{-1})\Lambda (x),\  g\in G\ .
        $$

        \noi We have proved:

        \begin{lemma} \label{model} We have a natural isomorphism of $G$-modules
          $$
          {\tilde \VV}\simeq \{\, \Lambda\in {\rm Ind}_J^G\, {\tilde \lambda}\ ; \ \varphi \star\Lambda
          =\chi (\varphi )\Lambda , \ \varphi \in \HH_\lambda \, \}
          $$
        \end{lemma}

        In order to produce an explicit matrix coefficient of $\pi$, we have to choose two vectors in
        $\VV$ and $\tilde \VV$ respectively. First fix a vector $w$ in the space $W$ of $\lambda$. Let
        $f_w \in \XX= {\rm ind}_J^G \, \lambda$ be the function with support $J$, given by $f_{w}(j)=
        \lambda (j).w$, $j\in J$. Its image ${\bar f}_w$ in $\XX /\XX_\chi$ is our favourite vector in
        $\VV$.

        \smallskip
        
        To produce vectors in ${\tilde \VV}$ we proceed as follows. Let ${\bar \HH}_\lambda$
        be the space of functions $\Psi$~: $G\lra {\rm End}_\CC ({\tilde W})$ satisfying
        $\Psi (j_1 gj_2 ) = {\tilde \lambda}(j_1 )\circ \Psi (g)\circ {\tilde \lambda}(j_2 )$, $j_1$,
        $j_2\in J$, $g\in G$ (with no condition on support).  Let ${\tilde w}\in {\tilde W}$ and
        $\Psi\in {\bar \HH}_\lambda$. Then the function $\Lambda_{\Psi, {\tilde w}}$~: $G\lra {\tilde W}$
        defined by $\Lambda_{\Psi ,{\tilde w}}(g) = \Psi (g)({\tilde w})$ is by construction
        an element of ${\rm Ind}_J^G \, {\tilde \lambda}$. 

        For $\varphi\in \HH_\lambda$, and $\Psi\in {\bar \HH}_\lambda$ define the convolution $\varphi\star
        \Psi$ by the {\it finite} sum
        $$
        (\varphi \star \Psi )(g) = \sum_{x\in J\backslash G} \varphi (gx^{-1})\Psi (x)\ , g\in G\ .
        $$
        
        \noi This is indeed an element of ${\bar \HH}_\lambda$.

        \begin{lemma}  \label{relation}
          Let $\Psi\in {\bar \HH}_\lambda$ be a function satisfying $\varphi \star \Psi =
          \chi (\varphi ) \Psi$, for all $\varphi \in \HH_\lambda$.
          Then, in the model of ${\tilde \VV}$ given in Lemma \ref{model}, we have
          $\Lambda_{\Psi , {\tilde w}} \in {\tilde \VV}$.
          \end{lemma}

        \noi {\it Proof}. If $\Psi\in {\bar \HH}_\lambda$ satisfies the assumption of the lemma, for
        $\varphi\in \HH_\lambda$ and $g\in G$, we have

    \begin{align*}
          (\varphi\star \Lambda_{\Psi, {\tilde w}})  (g) & = & \sum_{x\in J\backslash G}\varphi (gx^{-1})
          \Lambda_{\Psi, {\tilde w}} (x) \\
          &=&   \sum_{x\in J\backslash G}\varphi (gx^{-1}) \Psi (x)({\tilde w})\\
          &=&   ( \sum_{x\in J\backslash G}\varphi (gx^{-1}) \Psi (x)) ({\tilde w})\\
          &=& (\varphi\star \Psi )(g)({\tilde w})\\
          &=& \chi (\varphi )\Lambda_{\Psi ,{\tilde w}}(g)
        \end{align*}

     \noi as required.

     Let us write ${}_\chi{\bar \HH}_\lambda$ for the $\CC$-vector
     space given by
     $$
     \left\{ \Psi\in {\bar \HH}_\lambda \ ; \ \varphi \star \Psi =
     \chi (\varphi )\Psi , \ \varphi\in \HH_\lambda\right\}
     $$

     \noi One easily check that, for all $\Psi\in {}_\chi{\bar
       \HH}_\lambda$, the map $L_\psi$~: ${\tilde W} \ni w\mapsto
     \Lambda_{\Psi , w}\in {\tilde \VV}$ is $J$-equivariant. We
     therefore have a linear map $L$~: ${}_\chi {\bar \HH}_{\lambda}
     \lra {\rm Hom}_J \, ({\tilde W}, {\tilde \VV} )$. Moreover it is a
       straightforward consequence of the definition of $\Lambda_{\Psi
         ,w}$ that $L$ is injective.

       \begin{proposition}\label{mult1} We have ${\rm dim}_{\CC}\,
         {}_\chi {\bar \HH}_\lambda \leqslant 1$.
       \end{proposition}
 \noi {\it Proof}. The spaces ${\rm Hom}_J \, ({\tilde W}, {\tilde \VV} )$, ${\rm Hom}_J \, (\pi , \lambda )$ have the same dimension. But
   it is a standard fact (e.g. see \cite{BH} {\S}3, Lemma 7.(2))
   that the simple type $\lambda$ occurs with
   multiplicity $1$ in the discrete series representations $\pi$.
   \smallskip

   We shall see later that   ${}_\chi {\bar \HH}_\lambda$ is
   non-trivial.

     \begin{proposition} \label{formula1} Let $w\in \Weyl$, ${\tilde w}\in {\tilde W}$.
       Let $\Psi\in {\bar \HH}_\lambda$ satisfy
       $\varphi \star \Psi =\chi (\varphi )\Psi$, for all $\varphi\in \HH_\lambda$. Let
       $c_{\Psi ,w,{\tilde w}}$ be the matrix coefficient of $\pi$ attached to the vector ${\bar f}_w\in\VV$ and
       to the linear form $\Lambda_{\Psi ,{\tilde w}}\in {\tilde
         \VV}$.

       (a) We have the formula:
       $$
       c_{\Psi , w,{\tilde w}}(g)= \langle w, \Psi(g^{-1}) .{\tilde w}\rangle_{\tilde \lambda} , \ g\in G\ .
       $$

       (b) Assume $\Psi\not= 0$. There exist $J$-equivariant
       embeddings $W \subset \VV$ and ${\tilde W}\subset {\tilde \VV}$
       such that $c_{\Psi , w,{\tilde w}}$ is the matrix coefficient of $\pi$
       attached to the vector $w\in \VV$ and to the linear form
       ${\tilde w}\in {\tilde \VV}$. 
    \end{proposition}

     \noi {\it Proof}.  (a) Let $g\in G$. By definition we have
     \begin{align*}
       c_{\Psi, w, {\tilde w}} (g) & = & \sum_{x\in J\backslash G}\langle f_w (xg), \Lambda_{\Psi ,{\tilde w}}(x)
       \rangle_{\tilde \lambda} \\
       &=& \sum_{x\in J\backslash G} \langle f_{w}(xg), \Psi (x)({\tilde w})\rangle_{\tilde \lambda}
     \end{align*}

     \noi Since the support of $f_{w}$ is contained in $J$, only the coset $Jg^{-1}$ may contribute to the sum,
     whence we obtain:
     
 \begin{align*}
   c_{\Psi, w, {\tilde w}} (g) & = & \langle f_{w}(1), \Psi (g^{-1})({\tilde w})\rangle_{\tilde \lambda}\\
   &=& \langle w , \Psi (g^{-1})({\tilde w})\rangle_{\lambda}
 \end{align*}
 \noi as required.

\noi (b) Indeed $W\ni w\mapsto {\bar f}_w \in \VV$ and ${\tilde W}\ni
     {\tilde w}\mapsto \Lambda_{\Psi ,{\tilde w}}\in {\tilde \VV}$ are non-zero and
     $J$-equivariant.

 \medskip

     For $\Psi\in {\bar \HH}_\lambda$, the formula ${\tilde \Psi}(g)= {}^t \Psi (g^{-1})$, $g\in G$,
     defines an element of the $\CC$-vector space  ${\bar \HH}_{\tilde \lambda}$ of functions
     $\Psi'$~: $G\lra {\rm End}_\CC (W)$ satisfying  $\Psi' (j_1 g j_2 )=
     \lambda (j_1 )\circ \Psi '(g)\circ \lambda (j_2)$ for all $j_1 , \ j_2 \in J$, $g\in G$. 
     For $\Psi'\in {\bar \HH}_{\tilde \lambda}$ and $\varphi '\in \HH_{\tilde \lambda}$, we define
     the convolution $\Psi'\star \varphi$ by the {\it finite} sum
     $$
     \Psi' \star \varphi '(g)=\sum_{G/J} \Psi '(x)\varphi '(x^{-1}g)
     $$

     It is straightforward to check that an element $\Psi$ of ${\bar \HH}_{\lambda}$ satisfies
     $\varphi \star \Psi = \chi (\varphi )\Psi$, $\varphi \in \HH_\lambda$,
     if, and only if, ${\tilde \Psi}$ satisfies ${\tilde \Psi}\star \varphi ' ={\tilde \chi}(\varphi ')
     {\tilde \Psi}$, $\varphi '\in {\HH}_{\tilde \lambda}$, where ${\tilde \chi}$ is the
     character of $\HH_{\tilde \lambda}$ defined by ${\tilde \chi}({\tilde \varphi })=\chi (\varphi )$,
     $\varphi \in \HH_\lambda$.

     We can therefore reformulate Proposition \ref{formula1} as follows.

     \begin{proposition} \label{formula2} Let $w\in \Weyl$ and ${\tilde w}\in {\tilde W}$. Let $\Psi$
       be an element of ${\bar \HH}_{\tilde \lambda}$ satisfying $\Psi\star \varphi =
       {\tilde \chi}(\varphi )\Psi$, for all $\varphi \in \HH_{\tilde \lambda}$. Then the formula
       $$
       c_{\Psi ,w,{\tilde w}}(g) = \langle \Psi (g)(w), {\tilde w}\rangle_{\lambda}, \ g\in G ,
       $$
       defines a matrix coefficient of $\pi$.

       Moreover the support of that matrix coefficient is contained in
       $J\Weyl J$.
     \end{proposition}

     \noi {\it Proof}. Only the second assertion  deserves to be proved.
     The pair $(J,\lambda )$ is a {\it  simple type} in the sense of \cite{BK1}
     Definition (5.5.10). By {\it op. cit.} Proposition (5.5.11), if $g\in G$
     intertwine $(J,\lambda )$, that is if the intertwing space
     ${\rm Hom}_{J\cap {}^g J} (\lambda , {}^g\lambda )$ is not trivial, then $g$
     must lie in $J\Weyl J$. Moreover easily adapting the proof of {\it loc. cit.}
     Proposition (4.1.1),  if for $g\in G$, we have $\Psi (g)\not= 0$,
     then $g$ intertwines $(J,\lambda )$. As a consequence the support
     of $\Psi$, whence that of $c_{\Psi ,w,{\tilde w}}$,  is contained
     in $J\Weyl J$.  
     \medskip

     The rest of the section is devoted to constructing  a function $\Psi
     \in {\bar \HH}_{\tilde \lambda}$ satisfying the hypothesis of \ref{formula2}.

     Observe that  any element $\varphi$ of ${\bar \HH}_\lambda$
     may be written
     $$
     \varphi = \sum_{w\in \Weyl_0 , \ k\in \ZZ} c_{w, k}\,  \varphi_{\Pi}^k \star \varphi_{[w]}
     $$
     \noi for some complex numbers $c_{w, k}$, $w\in \Weyl_0$, $k\in \ZZ$, and where the {\it formal} sum
     $$
     \sum_{w\in \Weyl_0 , \ k\in \ZZ} c_{w, k}\,  \varphi_{\Pi}^k \star \varphi_{[w]}
     $$
     \noi denotes the element  of ${\bar \HH}_\lambda$,
     with support in $J\Weyl J$, and whose value at $\Pi^k w$, $k\in
     \ZZ$, $w\in \Weyl_0$, is $c_{w,k} (\varphi_\Pi^k\star \varphi_{[w]})
     (\Pi^k w)$.

     \begin{theorem} \label{TheFormula} Define an element of $\HH_\lambda$ by the formula
       $$
       \Psi_0 = \sum_{w\in \Weyl_0 , \ k\in \ZZ} (-\frac{1}{q_1})^{l(w)}\, \chi (\varphi_\Pi)^{-k}\, \varphi_\Pi^k  \star\varphi_{[w]}
       $$
       \noi where $q_1 =q^f$.  Then $\Psi_0$ satisfies the assumption  of Proposition \ref{formula1}:
       \begin{equation}
       \label{EqPsi}
       \varphi\star \Psi_0  =\chi (\varphi )\Psi_0 \ , \ {\rm for\  all}\ \varphi \in \HH_\lambda
       \end{equation}
     \end{theorem}

     \begin{corollary}
       We have ${\rm dim}_{\CC}\,   {}_\chi {\bar \HH}_\lambda  =1$.
      \end{corollary} 
     \noi {\it Proof}. Indeed by  Proposition \ref{mult1} we have ${\rm
       dim}_{\CC}\,   {}_\chi {\bar \HH}_\lambda  \leqslant 1$. On the
     other hand the previous theorem exhibits a non-zero element of
     ${\rm dim}_{\CC}\,   {}_\chi {\bar \HH}_\lambda  =1$.
     \medskip

     The proof of the theorem  consists in several steps.

     \begin{lemma} For all $\varphi\in \HH_\lambda$, we have
       $$
       \varphi\star \Psi_0  = 
\sum_{w\in \Weyl_0 , \ k\in \ZZ} (-\frac{1}{{q_1}})^{l(w)}\, \chi (\varphi_\Pi)^{-k}\, \varphi\star \varphi_\Pi^k  \star\varphi_{[w]}
$$
     \end{lemma}
     \noi {\it Proof}. We have to prove that the LHS and the RHS both
     evaluated at $g\in G$ coincide for all $g\in J\Weyl J$.
     But after evaluation at $g$, for any $g\in J\Weyl J$, for a matter of support, only a finite number of terms
     may contribute to the infinite formal sum,  so that our result follows from the usual bilinear property of the
     convolution product in $\HH_\lambda$.
     \medskip

     It suffices to prove Equality \ref{EqPsi} for $\varphi = \varphi_\Pi^l$, $l\in \ZZ$, and for
     $\varphi = \varphi_{[s_i]}$, $i=1,...,e-1$.

     We first deal with $\varphi$ of the form $\varphi_\Pi^l$, $l\in \ZZ$. Thanks to the last lemma, we have:
     $$
     \varphi_\Pi^l \star\Psi_0 =
     \sum_{w\in \Weyl_0 , \ k\in \ZZ} (-\frac{1}{{q_1}})^{l(w)}\, \chi (\varphi_\Pi)^{-k}\, \varphi_\Pi^{k+l}  \star\varphi_{[w]}
     $$
     \noi and our result follows by the change of variable $k' =k+l$.
     \medskip

     We now prove Equality \ref{EqPsi} for $\varphi =\varphi_{[s_i ]} $, $i=0,...,e-1$.
     First recall that $\chi (\varphi_{[s_i ]})=-1$
     so that we must prove that $\varphi_{[s_i ]}\star \Psi_0 =-\Psi_0$.

     We shall need:

     \begin{lemma} \label{Psi0} We have the equality
       $$
       \Psi_0 = \sum_{w\in \Weyl_0 , \ k\in \ZZ} (-\frac{1}{{q_1}})^{l(w)}\, \chi (\varphi_\Pi)^{-k}\, \varphi_{[w]}\star \varphi_\Pi^k 
       $$
     \end{lemma}
     
     \noi {\it Proof}.  First,  for $k\in \ZZ$, $w\in \Weyl_0$, we have
     $\varphi_\Pi^k \varphi_{[w]}\varphi_\Pi^{-k}= \varphi_{\Pi^k w\Pi^{-k}}$. This can be proved by induction on $l(w)$,
     reducing to the case $w\in \{ s_0 ,...,s_{e-1}\}$, where the equality follows from a straightforward calculation. 
     Now our result follows using the fact that $\Pi^k$ normalizes
     $\Weyl_0$ and by a suitable change of variable.
     \medskip
     
     Let $\Weyl_1$ be the set of $w\in \Weyl_0$ such that $l(s_i w )=l(w)+1$. Recall (Proposition \ref{HeckeCalc})
     that $w\in \Weyl_1$ if and only if
    for  any reduced expression $w=w_1 w_2 \cdots w_l$, $l=l(w)$, we have $w_1\not= s_i$. Let $\Weyl_2$ be the complement
     $\Weyl_0\backslash \Weyl_1$, so that we have $\Weyl_2 = s_i \Weyl_1$.  We may write
     \begin{align*}
       \varphi_{[s_i ]}\star \Psi_0 &=&
       \sum_{k\in \ZZ} \left( \sum_{w\in \Weyl_0} (-\frac{1}{{q_1}})^{l(w)}\, \chi (\varphi_\Pi)^{-k}\, \varphi_{[s_i ]}\star
       \varphi_{[w]} \right) \star \varphi_\Pi^k 
        \end{align*}

     \noi For $k\in \ZZ$, abbreviate:
     $$
     S_k = \sum_{w\in \Weyl_0} (-\frac{1}{{q_1}})^{l(w)}\, \chi (\varphi_\Pi)^{-k}\, \varphi_{[s_i ]}\star
       \varphi_{[w]} 
       $$
       \noi Using the multiplication rules of Proposition \ref{HeckeCalc2}, as well as
       suitable changes of variable, we successively write:
       \begin{flalign*}
         S_k &=  \sum_{w\in \Weyl_1} (-\frac{1}{{q_1}})^{l(w)}\, \chi (\varphi_\Pi)^{-k}\, \varphi_{[s_i ]}\star
         \varphi_{[w]}\\
         &
         +\sum_{w\in \Weyl_2} (-\frac{1}{{q_1}})^{l(w)}\, \chi (\varphi_\Pi)^{-k}\, \varphi_{[s_i ]}\star
       \varphi_{[w]} \\
       & = 
      \sum_{w\in \Weyl_1} (-\frac{1}{{q_1}})^{l(w)}\, \chi (\varphi_\Pi)^{-k}\, \star
      \varphi_{[s_i w]}\\
     & +\sum_{w'\in \Weyl_1} (-\frac{1}{{q_1}})^{l(w)}\, \chi (\varphi_\Pi)^{-k}\, \varphi_{[s_i ]} \star \varphi_{[s_i w']}\\
      & =
      \sum_{w\in \Weyl_2} (-\frac{1}{{q_1}})^{l(w)-1}\, \chi (\varphi_\Pi)^{-k}\, 
      \varphi_{[w]}\\
     & +\sum_{w'\in \Weyl_1} (-\frac{1}{{q_1}})^{l(s_i w')}\, \chi (\varphi_\Pi)^{-k} ({q_1}\varphi_{[w']}+({q_1}-1)\varphi_{[s_i w']})\\
      & =
      - {q_1}\sum_{w\in \Weyl_2} (-\frac{1}{{q_1}})^{l(w)}\, \chi (\varphi_\Pi)^{-k}\, \varphi_{[w]}\\
      & +{q_1}\sum_{w'\in \Weyl_1} (-\frac{1}{{q_1}})^{l(s_i w')}\, \chi (\varphi_\Pi)^{-k} \varphi_{[w']}\\
      & +({q_1}-1)
      \sum_{w'\in \Weyl_1}  (-\frac{1}{{q_1}})^{l(s_i w')}\, \chi (\varphi_\Pi)^{-k}\varphi_{[s_i w']}\\
      & = - {q_1}\sum_{w\in \Weyl_2} (-\frac{1}{{q_1}})^{l(w)}\, \chi (\varphi_\Pi)^{-k}\, \varphi_{[w]}\\
    &  -\sum_{w'\in \Weyl_1} (-\frac{1}{{q_1}})^{l( w')}\, \chi (\varphi_\Pi)^{-k} \varphi_{[w']}\\
    &  +({q_1}-1)\sum_{w\in \Weyl_2}  (-\frac{1}{{q_1}})^{l(w )}\, \chi (\varphi_\Pi)^{-k}\varphi_{[w ]}\\
      &=
      - \sum_{w\in \Weyl_0} (-\frac{1}{{q_1}})^{l(w)}\, \chi (\varphi_\Pi)^{-k}\,
       \varphi_{[w]} 
       \end{flalign*}
       and our result follows. This finishes the proof of Theorem \ref{TheFormula}.
       \medskip

       \noi {\it Remark}. Let us explain to the reader how we guessed the formula giving
       $\Psi_0$. Let ${\mathbb H}$ be a connected semisimple group defined and split  over $F$. Assume for simplicity sake that
       ${\mathbb H}$ is simply connected and that its relative root system is irreducible. Let $T$ be a maximal split torus
       of $T$, $\mathcal A$ be the apartement of the affine building of ${\mathbb H}$ attached to $T$. Fix a chamber $C$ of
       $\mathcal A$ and
       let $I$ be the Iwahori subgroup of ${\mathbb H} (F)$ fixing $C$.
       Let $\Weyl$ be the affine Weyl group of ${\mathbb H}$ attached to $T$.
       Recall that we have the decomposition ${\mathbb H} (F)=\ds \sqcup_{w\in \Weyl} IwI$. 
       The chamber $C$ determines a generating set of involutions $S$  of $\Weyl$. The pair $(\Weyl,S)$ is a Coxeter system;  let us
       denote by $l$ its length function. Let $\St_H$ be the Steinberg representation of $H={\mathbb H} (F)$. Finally let $c$ be the
       matrix coefficient of $\St_H$ attached to an Iwahori fixed vector of $\St_H$ and an Iwahori fixed linear form on $\St_H$.
       Then  $c$ is given by the following  formula (\cite{Bo} Equality (3) page 252, proof of Proposition 5.3):
       $$
       c(k_1 wk_2 ) = C.  (\frac{-1}{q})^{l(w)}, \ k_1 ,k_2 \in I, \ w\in \Weyl\ .
       $$
       \noi for some constant $C$. 

       \medskip

       Using Proposition \ref{formula2}, Theorem \ref{TheFormula} may be rewritten as follows.

       \begin{theorem}\label{TheFormula2} Define an element ${\tilde \Psi}_0$ of $\HH_{\tilde \lambda}$ by
       $$ 
         {\tilde \Psi_0} = \sum_{w\in \Weyl_0 , \ k\in \ZZ} (-\frac{1}{{q_1}})^{l(w)}\, \chi (\varphi_\Pi)^{-k}\,
         {\tilde \varphi}_{[w]} \star {\tilde \varphi}_\Pi^{-k}
         $$
         \noi Fiw $u\in W$, ${\tilde u}\in {\tilde W}$. Then the formula
         $$
         c_{u, {\tilde u}, {\tilde \Psi_0}}(g)=\langle {\tilde \Psi}_0 (g).u, {\tilde u}\rangle_{\lambda}, \ g\in G
         $$
         \noi defines a matrix coefficient of $\St_{e}(\tau)$. Its support is
         contained in $J\Weyl J$.
       \end{theorem}

        In fact most of the results of this section are not special to
        ${\rm GL}_n$ and extend, with the same proof, to much more
        general situations.
        \smallskip

         Assume that $\mathbb G$ is a connected reductive algebraic
         group defined over $F$. Fix a type $(J,\lambda )$ in
         $G={\mathbb G}(F)$ and,   with the notation of the present
         section, an irreducible representation $\pi$ in
         $S_\lambda (G)$. Let us make the following assumption:
         $$
         {\rm (H)}\quad M_\lambda (\pi ) = {\rm Hom}_J \, (\lambda
         ,\pi )\ {\rm is}\ 1{\rm -dimensional.}
         $$

         \noi Let $\chi$ be the character of $\HH (G,\lambda )$
         afforded by the module $M_\lambda (\pi )$. Define ${\bar
           \HH}_\lambda$ and ${}_\chi {\bar \HH}_\lambda$ as before. 

         \begin{theorem} (i) The space ${}_\chi {\bar \HH}_\lambda$ is
           of dimension $\leqslant 1$.

           \noi (ii) For any $\Psi\in {}_\chi {\bar \HH}_\lambda$,
           $w\in \lambda$ and ${\tilde w}\in {\tilde \lambda}$, the
           formula
           $$
           c_{\Psi ,w,{\tilde w}}(g) = \langle \Psi (g)(w),{\tilde
             w}\rangle_\lambda\ , g\in G
           $$
           \noi defines a matrix coefficient of $\pi$.
           \end{theorem}

          \noi {\it Remark}. If one moreover knows that $\HH
          (G,\lambda )$ is isomorphic to the Iwahori-Hecke algebra of
          some other reductive group defined over $F$, in a support
          preserving way, then one can exhibit a non zero element of
          ${}_\chi {\bar \HH}_\lambda$ using the method of this
          section. However we shall not push this further in this
          article.

       \section{Strongly discrete symmetric spaces and distinction}

       The notation in this section is independent of the previous ones.
       \smallskip

       Let $G/H$ be a $p$-adic reductive symmetric space: $G$ is the group of $F$-rational points of a  connected
       reductive group $\GG$ defined over $F$ and $H={\mathbb H}(F)$, where ${\mathbb H}$ is the reductive group
       $\GG^\theta$, for some $F$-rational involution $\theta$ of $\GG$.

       Let $\chi$ be a smooth abelian character of $H$.
       Recall that a smooth representation $(\pi , \VV )$ of $G$ is said to be {\it $\chi$-distinguished} if there
       is a non-zero linear form $l$~: $\VV\lra \CC$ such that
       $$
       l(\pi (g).v) = \chi (g)\, l(v), \ g\in G, \ v\in \VV .
       $$

       \noi If $\chi$ is the trivial character, one says that $\pi$ is {\it $H$-distinguished}.
       \medskip

       Let ${\mathbb A}_G$ be the maximal split torus
       in the center of $\GG$, ${\mathbb A}_G^+$ be the connected component of ${\mathbb A}_G \cap H$, and set
       $A_G^+ = {\mathbb A}_G^+ (F)$.  Fix Haar measures $d\mu_H$ and $d\mu_A$ on $H$ and $A_G^+$ and let $d{\mu}_{H/A_G^+}$ be the
       quotient Haar measure on $H/A_G^+$.
       \smallskip
       
       Let $(\pi , \VV)$ be an irreducible smooth representation of $G$ whose central character $\omega_\pi$
       is trivial on $A_G^+$. For $v\in \VV$, ${\tilde v}\in {\tilde \VV}$ one consider the matrix coefficient
       $c_{\pi ,v,{\tilde v}}$ defines by $c_{\pi, v,{\tilde v}}(g)=\langle \pi (g).v,{\tilde v}\rangle$.
       \smallskip
       
         Following Sakellaridis and Venkatesh \cite{SV}, one says that $G/H$ is a {\it strongly tempered}
         (resp. {\it strongly discrete}) symmetric space if for any tempered (resp. discrete series) irreducible
         representation $(\pi ,\VV )$ of $G$, whose central character is trivial on $A_G^+$, for all $v\in \VV$,
         for all ${\tilde v}\in {\tilde \VV}$, the matrix coefficient  $c_{\pi ,v,{\tilde v}}$ seen as a function on $H/A_G^+$ is
         $\mu_{H/A_G^+}$-integrable. Of course any strongly  tempered symmetric space is strongly discrete.
         \smallskip

         In the following sections, we shall use the following fact.

         \begin{theorem} \cite{GO} Let $F/F_0$ be a quadratic
           extension. Then the symmetric space ${\rm GL}_n (F)/{\rm GL}_n (F_0 )$ is strongly
           discrete.
        \end{theorem}

         We shall use this result to prove cases of distinctions. Indeed assume that $G/H$ is strongly discrete and let
         $\pi$ be an irreducible discrete series representation of $H$ satisfying $(\omega_\pi )_{\mid A_G^+}\equiv 1$. Fix
         ${\tilde v}$ in $\tilde \VV$. Then the formula
         $$
         l_{\tilde v}\ : \ v\in \VV\mapsto \int_{H/A_G^+} c_{\pi ,v,{\tilde v}} (g) \, d\mu_{H/A_G^+}({\dot g})
         $$
         \noi defines a $H$-invariant linear form on $\VV$, whence an element of ${\rm Hom}_H (\pi ,\CC )$. 
         So we have the following criterion:
         \medskip

        \begin{criterion} \label{Cri} If there exists a test vector for $l_{\tilde v}$, that is a vector $v\in \VV$ such that
         $$\int_{H/A_G^+} c_{\pi ,v,{\tilde v}} (g) \, d\mu_{H/A_G^+}({\dot g}) \not= 0
         $$
          \noi then $\pi$ is $H$-distinguished.
        \end{criterion}

\section{The Galois symmetric space ${\rm GL}_n (F)/{\rm GL}_n (F_0 )$}
We use the notation of the previous sections.

\smallskip
Let $F/F_0$ be a quadratic unramified extension. We denote by $\ofr_{F_0}$, $\pfr_{F_0}$, $k_{F_0}$, and $\varpi_{F_0}$ the
ring of integers of $F_0$, its maximal ideal, its residue field and the choice of a uniformizer respectively. 
Set $G={\rm GL}_n (F)$, $G_0 = {\rm GL}_n (F_0 )$ so that $G/G_0$ is
a reductive symmetric space. We denote by $k_{F_0}$ the residue field of $F_0$ and set $q_0 = {\#} k_{F_0}$, so that
$q=q_0^2$.  We set $J_0 =J\cap G_0$; this is the group
of units of the  principal hereditary order  $\Afr_0 := \Afr^\Gamma$ of ${\rm M}_n (F_0 )$.

We choose the uniformizer $\varpi_F$ so that $\varpi_F = \varpi_{F_0}\in F_0$. With that choice we have $\Pi\in G_0$ and
$\Weyl= \Pi^\ZZ \ltimes \Weyl_0 \subset G_0$.

\medskip

Irreducible discrete series representations of $G$ distinguished by $G_0$ are
classified (\cite{An},\cite{Ma}).

\begin{theorem} (\cite{An} Theorem 1.3,  \cite{Ma} Corollary 4.2)  Let $\eta_{F/F_0}$ be the quadratic character of $F_{0}^\times$ attached to the extension $F/F_0$ via
  local class field theory. See $\eta_{F/F_0}$ as a character of $G_0$ via $\eta_{F/F_0}(g)=\eta_{F/F_0}({\rm det}(g))$,
  $g\in G_0$. Let $\tau$ be an irreducible supercuspidal representation of ${\rm GL}_f (F_0 )$. Then $\St_e (\tau )$
  is distinguished by $G_0$ if and only if $\tau$ is $\eta_{F/F_0}^{e-1}$-distinguished.
\end{theorem}

The aim of this section is to give an entirely local proof of one implication of Matringe's theorem
when $\tau$ has level $0$. As a byproduct of the proof we shall
exhibit explicit test vectors.

\medskip

For simplicity sake
we only deal with the case $e$ odd. So we prove the following.

\begin{theorem} \label{DistGal} With the notation as above, assume that $\tau$ is of level $0$ and that $e$ is odd. Then
  if $\tau$ is ${\rm GL}_f (F_0 )$-distinguished, $\St_e (\tau )$ is $G_0$-distinguished.
\end{theorem}

\noi In this aim, we shall use Criterion \ref{Cri}.
\medskip

So let us fix an irreducible level $0$  supercuspidal representation $\tau$ of ${\rm GL}_f (F)$ as in {\S}2.
First we quote the following classical result without proof.

\begin{lemma} The representation $\tau$ is ${\rm GL}_f (F_0 )$-distinguished if and only if $(\omega_\tau )_{\mid F_0^\times}
  \equiv 1$ and $({\tilde \lambda_0}, X_0 )$ is ${\rm GL}_f (k_{F_0})$-distinguished.
\end{lemma}

\medskip

With the notation of the previous section, $A_G$ is the group $F^\times$ of scalars matrices in $G$, and $A_G^+ \simeq F_0^\times$
the group of scalar matrices in $G_0$. We fix a Haar measure $\mu_{G_0}$ on $G_0$ by imposing  that $\mu_{G_0}(J_0 )=1$, and
a Haar measure $\mu_{F_0^\times}$ on $F_0^\times$ normalized by imposing that the unit group $U(F_0 )$ has measure $1$. We let
$\mu_{G_0 /F_0^\times }$ denote the quotient measure on $G_0 /F_0^\times$. 
  \medskip
  
  Let us assume that $\tau$ is ${\rm GL}_f (F_0 )$-distinguished. Fix $v_0\in X_0$ (resp. ${\tilde v_0}\in {\tilde X}_0$),
  a non-zero vector fixed by ${\rm GL}_f(\ofr_{F_0})$. Set $v:= v_0 \otimes \cdots \otimes v_0 \in X$ and
  ${\tilde v}:= {\tilde v}_0 \otimes \cdots \otimes {\tilde v}_0\in {\tilde X}$. By construction both vectors are fixed
  by $J_0$; they are also fixed by the endomorphisms $t_i$, $i=0,...,e-1$, and $\Gamma$. Let $c_{v,{\tilde v}}$
  be the matrix coefficient of $\St_e (\tau )$ constructed in Theorem \ref{TheFormula2}:
  $$
  c_{v,{\tilde v}} (g) = \langle {\tilde \Psi}_0 (g).v, {\tilde v}\rangle_\lambda , \ g\in G
  $$

  We are going to establish the following formula.

  \begin{proposition} \label{FormulaInt} We have
    $$
    \int_{G_0 /F_0^\times} c_{v, {\tilde v}}(g)\, d\mu_{G_0 /F_0^\times}({\dot g})
    = e\, \langle v,{\tilde v}\rangle_\lambda \,  P_{\Weyl_0} (\frac{-1}{q_0^f})
    $$
    \noi where $P_{\Weyl_0}$ is the Poincar\'e series of the affine coxeter group
    $\Weyl_0$.
  \end{proposition}

  Recall that the Poincar\'e series $P_{\Weyl_0}$ is the formal series
  defined by
  $$
  P_{\Weyl_0}(X) = \sum_{w_0 \in \Weyl_0} X^{l(w_0 )} \ .
  $$

  \noi By \cite{Bott,St} it is given by the formula
  $$
  P_{\Weyl_0}(X) = \frac{1}{(1-X)^{e-1}} \prod_{i=1}^{e-1} \frac{1-X^{m_i}}{1-X^{m_i -1}}
  $$
  \noi where $m_1$, $m_2$,  ...,$m_{e-1}$ are the exponent of the spherical Coxeter group $S_e$
  (cf. \cite{Bou}, Chap. V, {\S}6, d\'efinition 2). Indeed we have $m_i = i$,
  $i=1,...,e-1$ (cf. \cite{Bou}, Planche I). In particular $P_{\Weyl_0}$ defines a non-vanishing function on the
  real open interval (-1,1). 
  \medskip

  Let us explain how Theorem \ref{DistGal} follows from the previous proposition. First,
  by the previous discussion, we have $P_{\Weyl_0}(-1/q_0^f )\not= 0$. We are thus reduced to proving
  that $\langle v , {\tilde v}\rangle_{\lambda}\not= 0$. First observe that
  $$
  \langle v , {\tilde v}\rangle_{\lambda} = \langle v, {\tilde v}\rangle_{{\tilde \lambda}_0^{e\otimes}}
  $$
  where ${\tilde \lambda}_0^{e\otimes}$ is the representation of ${\rm GL}_f (k_F )^{e\times}$ in $X$.
  Since ${\rm GL}_f (k_F )^{e\times}/{\rm GL}_f (k_{F_0})^{e\times}$ is a Gelfand pair,
  we have
  $$
  X^{{\rm GL}_f (k_{F_0})^{e\times}} = \CC\, v \  \ {\rm and} \ \ {\tilde X}^{{\rm GL}_f (k_{F_0})^{e\times}} = \CC\,{\tilde v}
  $$

  Now our claim follows from the following general result whose statement and proof were
  communicated to me by Dipendra Prasad.

  \begin{lemma} (D. Prasad) \label{dipendra} Let $H/K$ be a Gelfand pair where $H$ and $K$ are finite groups.
    Let $(\pi , \VV )$ be an irreducible representation of $G$ such that
    $\pi$ and $\tilde \pi$ are $H$-distinguished. Let $v$ and $\tilde v$
    be generators of the line $\VV^K$ and ${\tilde \VV}^K$ respectively.
    Then $\langle v, {\tilde v}\rangle\not= 0$.
   \end{lemma}

  {\it Proof}. The natural pairing $f$~: $V\times {\tilde V}\lra \CC$ is
  non-degenerate and $G$-invariant, whence $H$-invariant. Decompose $V$ as a direct
  sum of its isotypic components according to the action of $H$:
  $$
  \VV = \sum_{\rho \in {\hat H}} \VV_\rho
  $$
  \noi where ${\hat H}$ denotes the dual of $H$. Here $\VV_\rho$ is either trivial or isomorphic
  to $\rho$ as a $H$-module. Similarly we have $\ds {\tilde V}=\sum_{\rho\in {\hat H}} {\tilde \VV}_\rho$.
  Hence we have the decomposition:
  $$
  \VV\times {\tilde \VV} =\sum_{\rho , \tau\in {\hat H}} \VV_\rho\times {\tilde \VV}_\tau
  $$
  \noi Let ${\mathbf 1}_H$ be the trivial character of $H$. Then $f_{\mid \VV_{{\mathbf 1}_H} \times {\tilde \VV}_\tau}$
  is trivial as soon as $\tau\not\sim {\mathbf 1}_H$. Assume for a contradiction that $\langle v, {\tilde v}
  \rangle =0$. Then  $f_{\mid \VV_{{\mathbf 1}_H} \times \VV_{{\mathbf 1}_{H}}} =0$. This implies that the orthogonal
  of $v$ with respect to $f$ is ${\tilde \VV}$, a contradiction. 

  \medskip

  The rest of this section is devoted to proving Proposition \ref{FormulaInt}.
  \medskip

  Recall that the support of $c_{v, {\tilde v}}$ is contained in $J\Weyl J$. We first prove:

  \begin{lemma}
  \label{DoubleClass} We have $J\Weyl J\cap G_0 = J_0 \Weyl J_0$.
  \end{lemma}

   \noi {\it Proof}. We need the  following auxiliary lemma. 

   \begin{lemma} (\cite{Stev} Lemma 2.1) Let $\Gamma$ be a finite group acting on a group $H$ by group automorphisms.
     Set $K$ be a $\Gamma$-invariant subgroup of $H$ and set $H_0 := H^\Gamma$, $K_0 := K^\Gamma$. Finally let
     $w\in H_0$. Then if the first cohomology set $H^1 (\Gamma , K\cap wKw^{-1} )$ is trivial, we have
     $(KwK)\cap H_0 =  K_0 wK_0$.
   \end{lemma}

   \noi So we are reduced to proving that for all $w\in \Weyl$, we have
   $$
   H^1 ({\rm Gal}(F/F_0 ) , J\cap wJw^{-1})=\{ 1\} .
   $$
   
   \noi Fix $w\in \Weyl$. We have $\Afr = \Afr_0 \otimes_{\ofr_{F_0}}\ofr_F$. So $J\cap wJw^{-1}$ is the
   group of units of the $\ofr_F$-order $(\Afr_0 \cap w\Afr_0 w^{-1} )\otimes_{\ofr_{F_0}} \ofr_F$. Since $F/F_0$ is unramified,
   it is a classical fact  in Galois Cohomology that if $\mathfrak M$ is an $\ofr_{F_0}$-order in ${\rm M}_n (F_0 )$, then
   the following non-abelian cohomology set is trivial:
   $$
   H^1 ({\rm Gal}(F/F_0 ), ({\mathfrak M}\otimes_{\ofr_{F_0}} \ofr_F )^\times ) =\{ 1\}\ .
   $$
   \medskip
   
   \noi Let us write
   $$
   J_0 \Weyl J_0 =\bigsqcup_{k=0}^{e-1} \sum_{w_0 \in \Weyl_0} F_{0}^\times J_0 \Pi^k w_0 J_0
   $$

   We observe the function $c_{v, {\tilde v}}$ is
   bi-$J_0$-invariant. Moreover since $\tau$ is ${\rm GL}_f (F_0
   )$-distinguished, its central  character $\omega_\tau$ is trivial
   on $F_0^\times$. Since $\St_e (\tau )$ and $\tau$ share the same
   central character, it follows that any matrix coefficient of $\St_e (\tau
   )$ is $F_0^\times$-invariant.   As a consequence we may write:
   
   \begin{align*}
   \int_{G_0 /F_0^\times} c_{v,{\tilde v}} (g)\, d\mu_{G_0 /F_0^\times}({\dot g}) & = 
   \sum_{k=0}^{e-1}\sum_{w_0 \in \Weyl_0} \int_{J_0 w_0 \Pi^k  J_0} c_{v, {\tilde v}}(g)d\mu_{G_0} (g) \\
   & =\sum_{k=0}^{e-1} \sum_{w_0\in \Weyl_0} \mu_{G_0} (J_0 w_0 \Pi^k
   J_0 )c_{v,{\tilde v}}(w_0 \Pi^k  )
   \end{align*}

   \begin{lemma} \label{Measure}For all $w_0\in \Weyl_0$ and $k\in \ZZ$, we have
     $$
     \mu_{G_0} (J_0 w_0 \Pi^k  J_0 )= q_0^{f^2 l(w_0 )}\ .
     $$
   \end{lemma}
   \noi {\it Proof}. Since $\Pi$ normalizes $J_0$, we have $ \mu_{G_0}
   (J_0 w_0 \Pi^k  J_0 ) =  \mu_{G_0} (J_0 w_0 J_0 )$. Let $I_0$ be the
   standard Iwahori subgroup of ${\rm GL}_e (F_0 )$ formed of those
   matrices in ${\rm GL}_e ({\ofr}_{F_0})$ that are upper triangular
   modulo $\pfr_{F_0}$. If one considers $\Weyl_0$ as a subgroup of
   ${\rm GL}_e (F_0 )$ it coincides with the affine Weyl group of the
   diagonal torus. As a consequence, if $\mu_0$ denotes the Haar
   measure on ${\rm GL}_e (F_0 )$ normalized by $\mu_0 (I_0 )=1$, we
   have the Iwahori-Matsumoto formula $\mu_0 (I_0 wI_0 )=
   q_0^{l(w)}$, $w\in \Weyl_0$ 
   (\cite{IM} Proposition 3.2). Considering $J_0$ as a block version
   of the Iwahori subgroup $I_0$, the formula of the lemma is a block
   version of Iwahori and Matsumoto's formula, where one has to
   replace $q_0 =\vert \ofr_0 /\varpi_F \ofr_0\vert$ by $q_0^{f^2}
   = \vert {\rm M}_f (\ofr_0 )/\varpi_F {\rm M}_f (\ofr_0
   )\vert$. The proof is left to the reader. 
   \medskip
   
   We shall need a last lemma, whose proof we postpone to the end of
   that section.
   
   \begin{lemma} \label{PsiValue} For all $k\in \ZZ$ and $w_0\in \Weyl_0$, we have
     $$
     c_{v, {\tilde v}} (w_0 \Pi^k  )= (\frac{-1}{q_1})^{l(w_0 )} \, q^{-\frac{f(f-1)}{2}l(w_0 )} \langle v,{\tilde v}\rangle_{\lambda}\ .
     $$
   \end{lemma}
   \noi {\it Proof}.

   \medskip

   Thanks to the previous  series of lemmas, we now may  write:

   \begin{align*}
      \int_{G_0 /F_0^\times} c_{v,{\tilde v}} (g)\, d\mu_{G_0 /F_0^\times}({\dot g}) 
      & =\sum_{k=0}^{e-1} \sum_{w_0\in \Weyl_0} \mu_{G_0} (J_0 \Pi^k w_0 J_0 )\, c_{v,{\tilde v}}(\Pi^k w_0 )\\
      & = \sum_{k=0}^{e-1} \sum_{w_0\in \Weyl_0} q_0^{f^2 l(w_0 )}\,  (\frac{-1}{q_1})^{l(w_0 )} \,
        q^{-\frac{f(f-1)}{2}l(w_0 )} \, \langle v,{\tilde v}\rangle_{\lambda}\\
      & = e\, \sum_{w_0\in \Weyl_0} q_0^{f^2 l(w_0 )}\,  (\frac{-1}{q_0^{2f}})^{l(w_0 )} \,
          (q_0^2 )^{-\frac{f(f-1)}{2}l(w_0 )} \, \langle v,{\tilde v}\rangle_{\lambda}\\
          & = e\, \sum_{w_0\in \Weyl_0} (-1)^{l(w_0 )} q_{0}^{f^2 l(w_0 ) - 2fl(w_0 )-f^2 l(w_0 ) +fl(w_0 )} \, \langle v,{\tilde v}\rangle_{\lambda}\\
          & = e\, \sum_{w_0\in \Weyl_0} (\frac{-1}{q_0^f})^{l(w_0 )} \, \langle v,{\tilde v}\rangle_{\lambda} \\
          & = e\, P_{\Weyl_0}(\frac{-1}{q_0^f})\, \langle v, {\tilde v}\rangle_{\lambda}
   \end{align*}
   
   \noi as required. This finishes the proof of Proposition \ref{FormulaInt}. 
   \bigskip

\begin{theorem} \label{testvector} Assume that $e$ is odd. Fix any $J$-equivariant
  embedding $W\subset \VV = \St_e (\tau )$. Let $w\in W$ be a
  generator of $W^{J_0}$, where we recall that $J_0 = J\cap {\rm GL}_n
  (F_0 )$. Then for all non-zero linear form $\Phi\in {\rm Hom}_{{\rm GL}_n
    (F_0 )}\, (\St_e (\tau ), \CC )$, we have $\Phi (w)\not= 0$.
 \end{theorem} 

\noi {\it Proof}. It is a result of Flicker's \cite{F} that the pair
${\rm GL}_n (F) /{\rm GL}_n (F_0 )$ has the following {\it
  multiplicity $1$} property: for any irreducible smooth
representation $\sigma$ of ${\rm GL}_n (F)$, one has
$$
{\rm dim}_{\CC} \, {\rm Hom}_{{\rm GL}_n (F_0 )} (\sigma , \CC
)\leqslant 1\ .
$$

 \noi It follows that any  linear form $\Phi \in {\rm Hom}_{{\rm GL}_n
    (F_0 )}\, (\St_e (\tau ), \CC )$ is proportional to the form
 $\Phi_0$ given with the notation of {\S}4  by
 $$
 \Phi_0 (u) = \int_{G_0 /F_0^\times} c_{\pi , u,{\tilde v}} (g)\,
 d\mu_{G_0 /F_0^\times}({\dot g}), \ \ u\in \VV \ .
 $$

 \noi Now the fact that $w$ is a test vector for $\Phi_0$ follows from
 from Proposition \ref{formula1}(b) and from the formula of
 Proposition \ref{FormulaInt} together with Lemma \ref{dipendra}.
           
\bigskip

 \noi {\it Proof of Lemma \ref{PsiValue}}. By definition we have
 $c_{v,{\tilde v}}(w_0 \Pi^k )=\langle {\tilde \Psi}_0 (w_0 \Pi^k ).v ,{\tilde
   v}\rangle_\lambda$, with
 \begin{equation}\label{Psi1}
  {\tilde \Psi}_0 (w_0 \Pi^k ) = (-\frac{1}{q_1})^{l(w_0 )} \chi
  ({\tilde \varphi}_\Pi )^k \left( {\tilde \varphi}_{[w_0 ]} \star
  {\tilde \varphi}_\Pi^k \right) (w_0 \Pi^k )
 \end{equation}

 \noi By Theorem \ref{caractere}, we have $\chi ({\tilde \varphi}_\Pi ) =
 (-1)^{e-1} \omega_\tau ((-1)^{e-1} \varpi_{F_0})$. Since $e$ is odd
 and $\tau$ is ${\rm GL}_f (F_0 )$-distinguished, we obtain $\chi
 ({\tilde \varphi}_\Pi ) = 1$.

 A convolution calculation shows that ${\tilde \varphi}_{[w_0 ]} \star
 {\tilde \varphi}_\Pi^k (w_0 \Pi^k ) = {\tilde \varphi}_{[w_0 ]} (w_0
 ){\Gamma}^k$. Therefore from \eqref{Psi1}, we deduce:

 \begin{align}
   c_{v, {\tilde v} } (w_0 \Pi^k ) &  =  (-\frac{1}{q_1})^{l(w_0 )}\langle {\tilde \varphi}_{[w_0 ]} (w_0
     ){\Gamma}^k .v, {\tilde v}\rangle_\lambda\\
    \label{Coeff} & =  (-\frac{1}{q_1})^{l(w_0 )}\langle {\tilde \varphi}_{[w_0 ]} (w_0
     ). v, {\tilde v}\rangle_\lambda
 \end{align}

 \noi since $v$ is $\Gamma$-invariant. 
    
 Set $l= l(w_0 )$ and let $w_0 = s_{i_1} \cdots s_{i_l}$, $i_k\in \{0,
 ..., e-1\}$, $k=1,...,l$,  be a reduced expression of $w_0$.  We are going to prove that

 \begin{equation}\label{phi}
   {\tilde \varphi}_{[w_0 ]} (w_0 ) = \left(
   q^{-\frac{f(f-1)}{2}}\right)^{l(w_0 )} \, t_{i_1} \circ \cdots \circ
     t_{i_l}  
 \end{equation}

 \noi Then Lemma \ref{PsiValue} will follow since $v$ is
 $t_k$-invariant, for $k=0,...,e-1$. We proceed by induction on the
 length $l$ of $w_0$. If $l=1$, then by Lemma \ref{valeurshecke}, we
 have $\ds {\tilde \varphi}_{[w_0 ]} (w_0 ) = \omega_\tau (-1)
   q^{-\frac{f(f-1)}{2}} \, t_{i_1}$, with $\omega_\tau (-1) = 1$,
   since $\tau$ is ${\rm GL}_f (F_0 )$-distinguished, and the result
   follows. Assume that the result holds for elements of length $l$,
   for some $l\geqslant 1$, and let $w_0$ be of length $l+1$. Write
   $w_0 = s_i w_0 '$, with $w_0 '$ of length $l$ and $i\in \{ 0
   ,...,e-1\}$. It suffices to prove that

   \begin{equation}
     \label{induction}
   {\tilde \varphi}_{[w_0 ]}(w_0 ) =    q^{-\frac{f(f-1)}{2}} t_i \circ {\tilde
     \varphi}_{[w_0 ' ]}(w_0 ')
   \end{equation}

   By Cororollary \ref{HeckeCalc2}, we have ${\tilde \varphi}_{[w_0 ]} =
   {\tilde \varphi}_{[s_i ]} \star{\tilde \varphi}_{[w_0
       ']}$. Therefore:

   \begin{align}
{\tilde \varphi}_{[w_0 ]} (w_0 )& = \int_G  {\tilde \varphi}_{[s_i ]}
(w_0 g^{-1}) \circ {\tilde \varphi}_{[w_0 ']}(g)\, d\mu (g)\\
 & = \int_{Js_i Jw_0 \cap Jw_0 'J}  {\tilde \varphi}_{[s_i ]}
(w_0 g^{-1}) \circ {\tilde \varphi}_{[w_0 ']}(g)\, d\mu (g) \label{CalcInt}
   \end{align}

 \noi We have $Js_i Jw_0 \cap Jw_0 ' J=Js_i Js_i w_0 ' \cap Jw_0
 'J$. We have the containment
 $$
 Js_i Js_i \subset J\cup Js_i J
 $$
 \noi which is a ``block version'' of a classical axiom of Tits
 systems. We  deduce:

 \begin{align*}
Js_i Js_i w_0 & \subset Jw_0 ' \cup Js_i Jw_0 ' \\
   & \subset Jw_0 ' \cup Js_i Jw_0 'J \\
   & \subset Jw_0 ' \cup Js_i w_0 ' J
 \end{align*}

 \noi where, once again, the equality $Js_i Jw_0 'J = Js_i w_0 'J$
 follows from $l(s_i w_0 ')=l(w_0 ')+1$ and from a ``block version of
 a classical axiom of Tits systems (cf. \cite{BK1} Lemma (5.6.12) for
 a proof). Since $Js_i w_0 'J$ and $Jw_0 ' J$ are disjoint
 (cf. e.g. the proof of \cite{BK1} Proposition (5.5.16)), we deduce
 that $Js_i Jw_0  \cap Jw_0 ' J \subset Jw_0 ' = Js_i w_0$, and that  $Js_i
 Jw_0  \cap Jw_0 '  = Js_i w_0$, since the containment $Jw_0 ' \subset
 Js_i Jw_0  \cap Jw_0 '$ is obvious.

 So from  \ref{CalcInt}, we may  write:
 \begin{align*}
 {\tilde \varphi}_{[w_0 ]} (w_0 ) & = \int_{Js_i w_0 '}  {\tilde \varphi}_{[s_i ]}
 (w_0 g^{-1}) \circ {\tilde \varphi}_{[w_0 ']}(g)\, d\mu (g) \\
 & = {\tilde \varphi}_{[s_i ]} 
 (s_i ) \circ {\tilde \varphi}_{[w_0 ']}( w_0 ')
 \end{align*}
 This finishes the proof of \ref{induction} and of  Lemma \ref{PsiValue}.

Paul Broussous
\smallskip

Universit\'e de Poitiers 

Laboratoire de Math\'ematiques et Applications, UMR 7348

Site du Futuroscope - T\'el\'eport 2

11 Boulevard Marie et Pierre Curie

B\^atiment H3 - TSA 61125

86073 POITIERS CEDEX 9 

France
\smallskip

paul.broussous{@}math.univ-poitiers.fr

\end{document}